\documentclass[reqno]{amsart}
\usepackage{amsfonts}
\usepackage{amsmath,amsthm,amsfonts,amssymb}
\newtheorem{lemma}{Lemma}
\newtheorem{theorem}{Theorem}
\newtheorem{corollary}{Corollary}

\def\bl{\begin{lemma}}
\def\bt{\begin{theorem}}
\def\el{\end{lemma}}
\def\et{\end{theorem}}
\def\bp{\begin{proof}}
\def\ep{\end{proof}}
\def\bc{\begin{corollary}}
\def\ec{\end{corollary}}
\def\iy{\infty}
\def\te{\theta}

\def\mb{\mathbb}
\def\su{\subset}
\def\l{\lambda}

\def\o{\omega}
\def\O{\Omega}
\def\a{\alpha}
\def\b{\beta}
\def\p{\partial}
\def\su{\subset}

\def\s{\sigma}

\def\D{\Delta}
\def\g{\gamma}

\def\z{\zeta}

\def\-{\setminus}
\def\s{\sigma}

\def\vp{\varphi}
\def\ov{\overline}
\def\lt{\left}
\def\rt{\right}
\def\+{\bigcup}
\def\.{\bigcap}

\def\ll{\langle}
\def\rl{\rangle}

\title[The Schwarz-Pick lemma]
{The Schwarz-Pick lemma of high order\\ in several variables}
\thanks{Research supported by NSFC(China): 10171047, 10671093.}

\author {Shaoyu Dai, Huaihui Chen and Yifei Pan}

\address{Department of Mathematics, Nanjing Normal
University, Nanjing 210097, P.R.China}
\address{\it E-mail address: dymdsy@163.com}

\address{Department of Mathematics, Nanjing Normal
University, Nanjing 210097, P.R.China}
\address{\it E-mail address: hhchen@njnu.edu.cn}

\address{School of sciences,
NanChang University, Nanchang 330022, P.R.China}
\address{Department of Mathematical Sciences, Indiana University -
Purdue University Fort Wayne, Fort Wayne, IN 46805-1499}
\address{\it E-mail address: pan@ipfw.edu}

\begin{document}

\begin{abstract} We prove a high order Schwarz-Pick lemma for
mappings between unit balls in complex spaces in terms of the
Bergman metric. From this lemma, Schwarz-Pick estimates for partial
derivatives of arbitrary order of mappings are deduced.
\end{abstract}

\maketitle

\smallskip \noindent {\bf AMS Mathematics Subject Classifications
(2000) 32A10, 32F45.}


\section{Introduction}

Let $\mathbb{B}_n$ be the unit ball in the complex space
$\mathbb{C}^n$ of dimension $n$. The unit disk in the complex plane
is denoted by $\mathbb{D}$. For $z=(z_1,\cdots,z_n)$ and
$z'=(z'_1,\cdots,z'_n)\in\mb C^n$, denote $\ll z,z'\rl=z_1\ov
z_1'+\cdots+z_n\ov z_n'$ and $|z|=\ll z,z\rl^{1/2}$.

A multi-index $\a=(\a_{1},\cdots,\a_{n})$ of dimension $n$ consists
of $n$ non-negative integers $\a_{j}$, $1\leq j \leq n$, the degree
of a multi-index $\a$ is the sum
$|\a|=\overset{n}{\underset{j=1}{\sum}}\a_{j}$, and we denote
$\a!=\a_1!\cdots \a_n!$. For $z=(z_{1},\cdots,z_{n})\in
{\mathbb{C}}^{n}$ and a multi-index $\a=(\a_{1},\cdots,\a_{n})$, let
$z^\a=\overset{n}{\underset{j=1}{\prod}}{z_{j}}^{\a_{j}}$. A
holomorphic function $f$ on $\mb B_n$ can be expressed by
$f(z)=\sum\limits_\a a_\a z^\a$. For two multi-indexes
$\a=(\a_{1},\cdots,\a_{n})$ and $v=(v_{1},\cdots,v_{n})$, let
$v^\a=v_1^{\a_1},\cdots,v_n^{\a_n}$. Note that $v_j^{\a_j}=1$ if
$v_j=\a_j=0$.  Let $\O_{n,m}$ be the class of all holomorphic
mappings $f$ from $\mb{B}_n$ into $\mb{B}_m$. For $f\in\O_{n,m}$, if
$f=(f_1,\cdots,f_m)$, $f_j(z)=\sum\limits_{\a}a_{j,\a}z^{\a}$ for
$j=1,\cdots,m$, we denote $f(z)=\sum\limits_{\a}a_{\a}z^{\a}$, where
$a_\a=(a_{1,\a},\cdots,a_{m,\a})$.

For $f\in\O_{1,1}$, the classical Schwarz-Pick lemma says that
$$\frac{|f'(z)|}{1-|f(z)|^2}\le\frac1{1-|z|^2}$$
holds for $z\in\mb D$. Recently, the above inequality has been
generalized to the derivatives of arbitrary order by some authors
\cite{1,2,3}. The best result was proved in \cite{3}. It was proved
that
$$\frac{|f^{(k)}(z)|}{1-|f(z)|^2}\le(1+|z|)^{k-1}\cdot\frac{k!}
{(1-|z|^2)^k}\eqno(1.1)$$ holds for $f\in\O_{1,1}$, $k\ge1$ and
$z\in\mb D$. The equality in (1.1) may be attained if $z=0$, and the
equality statement has been established. If $k>1$ and $z\ne0$, (1.1)
is a strict inequality.

 Chen and Liu \cite{4} generalized (1.1) by proving
the following Schwarz-Pick estimate for partial derivatives of
arbitrary order of a function $f\in\O_{n,1}$:
$$\lt|\frac{\p^{|v|}f(z)}{\p{z_{1}}^{v_{1}}
...\p{z_{n}}^{v_{n}}}\rt|\leq n^{\frac{|v|}{2}}|v|!{n+|v|-1 \choose
n-1}^{n+2} \frac{1-|f(z)|^{2}}{(1-|z|^{2})^{|v|}}
(1+|z|)^{|v|-1}\eqno(1.2)$$ holds for any $z\in\mb B_n$ and
multi-index $v=(v_1,\cdots,v_n)\ne0$.

On the unit ball $\mb B_n$, the Bergman metric $H_z(\b,\b)$ may be
defined by
$$H_z(\b,\b)=\frac{(1-|z|^2)|\b|^2+|\ll \b,z\rl|^2}{(1-|z|^2)^2}
\quad\mbox{for}\ \ z\in\mb B_n,\ \b\in\mb C^n.$$ Commonly, there is
a factor $(n+1)/2$ in the definition of the Bergman metric. In spite
of ambiguity, we use the same notation for Bergman metrics in unit
balls of different dimensions. This metric is invariant under the
automorphism group of $\mb B_n$. For $f\in\O_{n,m}$, the
Schwarz-Pick lemma is formulated in terms of the Bergman metric (see
\cite{5}):
$$H_{f(z)}(f'(z)\b,f'(z)\b)\le H_z(\b,\b)\quad\mbox{for}\ \
z\in\mb B_n,\ \b\in\mb C^n.\eqno(1.3)$$ Here, $f'(z)$ is the
Jacobian matrix  of the  mapping $f$ at the point $z$, i.e.,
$f'(z)=\lt(\p f_j(z)/\p z_k\rt)_{1\le j\le m,1\le k\le n}$, and we
identify a point in complex space with a column matrix (column
vector) so that $f'(z)\b$ is the product of two matrixes. (1.3) is
precise, the equality holds for mappings in the automorphism group
of $\mb B_n$ if $m=n$.

The purpose of this paper is to generalize (1.3) to the high order
Fr\'echet derivatives of mappings in $\O_{n,m}$ as was done in
\cite{3} for the classical Schwarz-Pick lemma. For $f\in\O_{n,m}$,
$k\ge1$, $z\in\mb B_n$, the Fr\'echet derivative of $f$ at $z$ of
order $k$ is defined by
$$D_k(f,z,\b)=\sum_{|\a|=k}\frac{k!}{\a!}\frac{\p^kf(z)}
{\p z_1^{\a_1}\cdots\p z_n^{\a_n}}\b^\a,$$ where $\b\in\mb C_n$.
$D_k(f,z,1)=f^{(k)}(z)$ when $n=m=1$. With this notation, our main
result is expressed as follows:
\medskip

\noindent {\bf Theorem}:\quad {\it Let $f\in\O_{n,m}$. Then, for
$k\ge1$, $z\in\mb B_n$ and $\b\in\mb C^n\-\{0\}$, we have
\begin{align*}
H_{f(z)}(D_k&(f,z,\b),\ D_k(f,z,\b))\\
&\le k!^2\lt(1+\frac{|\ll\b,z\rl|}
{((1-|z|^2)|\b|^2+|\ll\b,z\rl|^2)^{1/2}}\rt)^{2(k-1)}(H_z(\b,\b))^k.
\tag{1.4}
\end{align*}
}
\medskip

\noindent (1.4) coincides with (1.1) or (1.3) if $n=m=1$ or $k=1$
respectively. Note that the factor preceding $(H_z(\b,\b))^k$ is
increasing with $|\ll\b,z\rl|$ from 0 to $|z|$.

As a consequence, we deduce from (1.4) a Schwarz-Pick estimate for
partial derivatives of a mapping $f\in\O_{n,m}$:
$$\lt|\lt\ll\frac{\p^{|v|}f(z)}{\p
z_1^{v_1}\cdots \p z_n^{v_n} },\ f(z)\rt\rl\rt|^2+(1-|f(z)|^2)\lt|
\frac{\p^{|v|}f(z)}{\p z_1^{v_1}\cdots \p z_n^{v_n}}\rt|^2$$
$$\le\frac{|v|^{|v|}}{v^v}\lt[v!(1+|z|)^{|v|-1}\cdot
\frac{1-|f(z)|^2}{(1-|z|^2)^{|v|}}\rt]^2\eqno(1.5)$$ holds for any
multi-index $v=(v_1,\cdots,v_n)\ne0$ and $z\in\mb B_n$. In
particular, if $f\in\O_{n,1}$, we have
$$\lt|\frac{\p^{|v|}f(z)}{\p z_1^{v_1}\cdots \p z_n^{v_n}}\rt|
\le\sqrt{\frac{|v|^{|v|}}{v^v}}v!(1+|z|)^{|v|-1}\cdot
\frac{1-|f(z)|^2}{(1-|z|^2)^{|v|}}.\eqno(1.6)$$  The equalities in
(1.5) and (1.6) may be attained if $z=0$ and the equality statement
is given. (1.6) is much better than (1.2) since the factor ${n+|v|-1
\choose n-1}^{n+2}$ is canceled,  $v!\le|v|!$ and
$\sqrt{|v|^{|v|}/v^v}\le n^{|v|/2}$ (the equality holds if and only
if $v_1=\cdots=v_n$).

For radial and normal partial derivatives, we have estimates  more
precise than (1.5) and (1.6). For $f\in\O_{n,1}$, we prove that
$$\lt|\frac{\p^{|v|}f(z)}{\p z_1^{v_1}\cdots \p z_n^{v_n}}\rt|
\le\sqrt{\frac{|v|^{|v|}}{v^v}}v!\mu(z)\cdot
\frac{1-|f(z)|^2}{(1-|z|^2)^{(v_1+|v|)/2}}\eqno(1.7)$$ holds for any
multi-index $v=(v_1,\cdots,v_n)\ne0$ and $z=(z_1,0,\cdots,0)\in\mb
B_n$, where $\mu(z)=(1+|z|)^{|v|-1}$ if $v_1=|v|$, and $\mu(z)$ is
the sum of terms $c_j|z|^j$ with $j\le v_1$ in $(1+|z|)^{|v|-1}$.

\section{Some lemmas}

The following results are known \cite{6}. For a point $a$ in a unit
ball, let
$$\vp_a(z)=\frac{a-P_az
-\sqrt{1-|a|^{2}}Q_az}{1-\langle z,a\rangle},$$ where $P_az=\ll
z,a\rl a/\ll a,a\rl$, $Q_az=z-P_az$. Note that $P_0(z)=0$. Then,
$\vp_a$ is  injective and maps the unit ball onto itself,
$$\vp_a(0)=a,\quad \vp_a(a)=0,\quad \vp_a=\vp_a^{-1},$$
and
$$\vp'_a(0)=-(1-|a|^2)P_a-(1-|a|^2)^{1/2}Q_a,$$
$$\vp'_a(a)=-\frac1{1-|a|^2}P_a-\frac1{(1-|a|^2)^{1/2}}Q_a.$$

\bl If $f(z)=\sum\limits_{\a}a_{\a}z^{\a}\in\O_{n,m}$, then
$$\sum_\a|a_\a|^2\lt|\b^{2\a}\rt|\leq1\eqno(2.1)$$
holds for $\b\in\p\mb B_n$. Further,
$$\sum_\a|a_\a|^2\cdot\frac{v^\a}{|v|^{|\a|}}\leq1\eqno(2.2)$$
holds for any multi-index $v=(v_1,\cdots,v_n)\ne0$. As a
consequence, we have
$$|a_v|\le\sqrt{\frac{|v|^{|v|}}{v^v}}.\eqno(2.3)$$
Further, if $v_j\ne0$ for $j=1,\cdots,n$, then the equality in (2.3)
holds only if $a_\a=0$ for $\a\ne v$. \el

\bp Let $\b=(\b_1,\cdots,\b_n)\in\p\mb B_n$ be fixed. For $0<\s<1$,
we have
\begin{equation*}
\begin{split}
1&\ge\frac1{(2\pi)^n}\int_0^{2\pi}\!\!\!\!\cdots\!\int_0^{2\pi}
|f(\s\b_1e^{i\te_1},\cdots,\s\b_ne^{i\te_n})|^2d\te_1\cdots d\te_n\\
&=\frac1{(2\pi)^n}\sum_{j=1}^m\int_0^{2\pi}\!\!\!\!\cdots\!\int_0^{2\pi}
|f_j(\s\b_1e^{i\te_1},\cdots,\s\b_ne^{i\te_n})|^2d\te_1\cdots
d\te_n\\
&=\sum_{j=1}^m\sum_\a|a_{j,\a}|^2\s^{2|\a|}|\b_1|^{2\a_1}\cdots
|\b_n|^{2\a_n}=\sum_\a|a_\a|^2\s^{2|\a|}|\b_1|^{2\a_1}\cdots
|\b_n|^{2\a_n}.
\end{split}
\end{equation*}
Letting $\s\to1$ gives (2.1). Thus, for given
$v=(v_1,\cdots,v_n)\ne0$, letting $\b_j=\sqrt{v_j/|v|}$ for
$j=1,\cdots,n$ in (2.1), we obtain (2.2). The lemma is proved. \ep

In the above proof, in order to get the best estimate (2.3) for
$a_v$, we deduce (2.2) by choosing $\b_j=\sqrt{v_j/|v|}$ in (2.1),
since the maximum $\underset{\b\in\p\mb
B_n}{\max}|\b^v|=\sqrt{\frac{v^v}{|v|^{|v|}}}$ is attained when
$\b_j=\sqrt{v_j/|v|}$ for $j=1,\cdots,n$.

\bl If  $f(z)=\sum\limits_{\a}a_{\a}z^{\a}\in\O_{n,m}$, then
$$\sum_{k=0}^\iy\lt|\sum_{|\a|=k}a_\a\b^a\rt|^2\le1\eqno(2.4)$$
holds for $\b\in\p \mb B_n$. \el

\bp For $\b\in\p \mb B_n$, let
$$h(\l)=f(\b\l)=\sum_{k=0}^\iy\lt(\sum_{|\a|=k}a_\a\b^\a\rt)\l^k,\quad
\l\in\mb D.$$ Then, $h(\mb D)\in\mb B_m$. Using (2.1), we obtain
(2.4). The lemma is proved. \ep

\bl Let $k\ge2$ be a positive integer and
$f(z)=\vp_{a}(bz^k)+g(z)$ for $z\in\mb D$, where $a\in\mb B_m$,
$b\in\p \mb B_m$ and
$$g(z)=\sum_{j=1}^{k-1}\sum_{n=0}^\infty a_{nk+j}z^{nk+j}$$
is a holomorphic mapping of $\mb D$ into $\mb C^m$. If $|f(z)|<1$
for $z\in\mb D$, then $g(z)\equiv0$.
\el

\bp Since $|g(z)|<1+|\vp_{a}(bz^k)|<2$, by Lemma 1, we have
$$\sum_{j=1}^{k-1}\sum_{n=0}^\infty|a_{nk+j}|^2\le4.$$
Thus, for $j=1,2,\cdots,k-1$, every component of the mapping
$$g_j(z)=\sum_{n=0}^\infty a_{nk+j}z^{nk+j}$$
is in the Hardy class $H^2$ and, consequently, for almost every
$\zeta\in\partial\mathbb D$, the radial limit
$\lim\limits_{z\to\zeta}g_j(z)$ exists for all $j$. Let $\zeta$ be
such a point and $\l=\vp_{a}(b\zeta^k)$. Obviously, $\l\in\p\mb
B_m$. Denote $\omega=e^{2\pi i/k}$. For $l=1,\cdots,k$, we have
$$\lim_{z\to\zeta}f(\omega^lz)
=\vp_{a}(b\zeta^k)+\sum_{j=1}^{k-1}\lim_{z\to\zeta}g_j(\omega^lz)
=\l+\sum_{j=1}^{k-1}\omega^{lj}\lim_{z\to\zeta}g_j(z),$$ and, since
$f(\mathbb D)\subset\mathbb B_m$,
$$\left|1+\sum_{j=1}^{k-1}\omega^{lj}\langle\lim_{z\to\zeta}
g_j(z),\l\rangle\right|\le\left|\lim_{z\to\zeta}f(\omega^lz)\right|\le1.$$

For $l=1,\cdots,k$, let
$$A_l=\sum_{j=1}^{k-1}\omega^{lj}\langle\lim_{z\to\zeta}g_j(z),\l\rangle.$$
Then, $|1+A_l|\le1$ and, consequently, $\mbox{Re}A_l\le0$ for
$l=1,\cdots,k$. However,
$$\sum_{l=1}^{k-1}A_l=\sum_{l=1}^{k-1}\sum_{j=1}^{k-1}\omega^{lj}\langle
\lim_{z\to\zeta}g_j(z),\l\rangle$$ $$=\sum_{j=1}^{k-1}\lt(\langle
\lim_{z\to\zeta}g_j(z),\l\rangle\sum_{l=1}^{k-1}\omega^{lj}\rt)
=-\sum_{j=1}^{k-1}\langle \lim_{z\to\zeta}g_j(z),\l\rangle=-A_k.$$
Thus, $\mbox{Re}A_k=0$. Noting that $|1+A_k|\le1$ we conclude that
$A_k=0$, i.e.,
$$\sum_{j=1}^{k-1}\langle\lim_{z\to\zeta}g_j(z),\l\rangle=0.$$
Thus,
$$1\ge\lt|\lim_{z\to\zeta}f(z)\rt|^2=
\lt|\l+\sum_{j=1}^{k-1}\lim_{z\to\zeta}g_j(z)\rt|^2
=1+\lt|\sum_{j=1}^{k-1}\lim_{z\to\zeta}g_j(z)\rt|^2
=1+\lt|\lim_{z\to\zeta}g(z)\rt|^2.$$ This shows that the radial
limit of every component of $g(z)$ is equal to $0$ at almost every
$\z\in\p\mb D$. According the general theory of $H^p$ spaces, we
conclude that $g(z)\equiv0$. The lemma is proved. \ep

\section{The partial derivatives at the origin}

\bt Let $f(z)=\sum\limits_{\a}a_{\a}z^{\a}\in\O_{n,m}$. Then,
$$\lt|\lt\ll\sum_{|\a|=k}a_\a\b^\a,\ a_0\rt\rl\rt|^2
+(1-|a_0|^2)\lt|\sum_{|\a|=k}a_\a\b^\a\rt|^2\le(1-|a_0|^2)^2.
\eqno(3.1)$$ holds for $k\ge1$ and $\b\in\p\mb B_n$. \et

\bp Let $k\ge1$ and $\b\in\p\mb B_n$ be given. If $a_0=0$, (3.1) is
a consequence of (2.4). Now, assume that $a_0\ne0$. Let
$$h(z)=\frac{1}{k}\sum_{l=1}^k f(e^{2l\pi i/k}z).$$
Then, $h(z)\in\O_{n,m}$, $h(0)=a_0$, and
$$h(z)=a_0+\sum_{m=1}^\iy\sum_{|\a|=mk}a_{\a}z^{\a}.$$
Let $\phi=\vp_{a_0}\circ h$. Obviously, $\phi\in\O_{n,m}$ and
$\phi(0)=0$. We have
\begin{equation*}
\begin{split}
\phi(z)&=\frac1{1-\langle
h(z),a_0\rangle}\lt(-(a_0/|a_0|^2)\sum_{m=1}^\iy\sum_{|\a|=mk}
\langle a_{\a},a_0\rangle z^{\a}\rt.\\
&\ \ \
\lt.-\sqrt{1-|a_0|^{2}}\sum_{m=1}^\iy\sum_{|\a|=mk}a_{\a}z^{\a}
+\sqrt{1-|a_0|^{2}}(a_0/|a_0|^2)\sum_{m=1}^\iy\sum_{|\a|=mk}
\langle a_{\a},a_0\rangle z^{\a}\rt)\\
&=-\frac1{1-\langle
h(z),a_0\rangle}\sum_{m=1}^\iy\sum_{|\a|=mk}\lt(\frac{\langle
a_{\a},a_0\rangle a_0}{1+\sqrt{1-|a_0|^{2}}}
+\sqrt{1-|a_0|^{2}}a_{\a}\rt)z^{\a}\\
&=-\frac1{1-|a_0|^2}\sum_{|\a|=k}\lt(\frac{\ll a_\a,a_0\rl
a_0}{1+\sqrt{1-|a_0|^{2}}} +\sqrt{1-|a_0|^{2}}a_\a\rt)z^{\a}
+\sum_{m=2}^\iy\sum_{|\a|=mk}c_\a z^{\a}.
\end{split}
\end{equation*}
Thus, using (2.4), we obtain
$$\frac1{(1-|a_0|^2)^2}\lt|\sum_{|\a|=k}\lt(\frac{\ll
a_{\a}\b^\a,a_0\rl a_0}{1+\sqrt{1-|a_0|^{2}}}
+\sqrt{1-|a_0|^{2}}a_{\a}\b^\a\rt)\rt|^2\le1.$$ A simple calculation
gives
\begin{equation*}
\begin{split}
&\lt|\sum_{|\a|=k}\lt(\frac{\ll a_{\a}\b^\a,a_0\rl
a_0}{1+\sqrt{1-|a_0|^{2}}}
+\sqrt{1-|a_0|^{2}}a_{\a}\b^\a\rt)\rt|^2\\
&\quad\quad\quad\quad =\lt|\frac{\lt\ll \sum_{|\a|=k}a_{\a}\b^\a,\
a_0\rt\rl a_0}{1+\sqrt{1-|a_0|^{2}}}
+\sqrt{1-|a_0|^{2}}\sum_{|\a|=k}a_{\a}\b^\a\rt|^2
\end{split}
\end{equation*}
\begin{equation*}
\begin{split}
&=\frac{\lt|\lt\ll\sum_{|\a|=k} a_\a\b^\a,\
a_0\rt\rl\rt|^2|a_0|^2}{\lt(1+\sqrt{1-|a_0|^2}\rt)^2}
+(1-|a_0|^2)\lt|\sum_{|\a|=k}a_{\a}\b^\a\rt|^2\\
&\qquad\qquad\qquad\qquad\qquad+\frac{2\sqrt{1-|a_0|^2}\lt|\lt\ll\sum_{|\a|=k}
a_\a\b^\a,
a_0\rt\rl\rt|^2}{1+\sqrt{1-|a_0|^2}}\\
&=\lt|\lt\ll\sum_{|\a|=k}a_\a\b^\a,\
a_0\rt\rl\rt|^2+(1-|a_0|^2)\lt|\sum_{|\a|=k}a_{\a}\b^\a\rt|^2.
\end{split}
\end{equation*}
This shows (3.1). The theorem is proved.
\ep

\bt Let $f(z)=\sum\limits_{\a}a_{\a}z^{\a}\in\O_{n,m}$. Then,
$$\lt|\lt\ll a_v,\ a_0\rt\rl\rt|^2
+(1-|a_0|^2)\lt|a_v\rt|^2\le\frac{|v|^{|v|}}{v^v}(1-|a_0|^2)^2.
\eqno(3.2)$$ holds for any multi-index $v\ne0$. Further, if the
equality holds for some $v=(v_1,\cdots,v_n)$ with $v_j\neq0$ for
$j=1,\cdots,n$, then
$$f(z)=a_0+\frac{a_v z^v}{1+\frac{\langle a_v, a_0\rangle
z^v}{1-|a_0|^2}}=a_0+a_v z^v+\cdots.\eqno(3.3)$$ Conversely, if
$v\neq0$, $a_0\in\mb B_m$ and $a_v\in\mb C^m$ satisfy the equality
in (3.2), then the mapping $f$ expressed by (3.3) belongs to
$\O_{n,m}$. \et

\bp Let $v=(v_1,\cdots,v_n)\ne0$ be given and $k=|v|$.  As in the
proof of the above theorem, consider $h$ and $\phi$. Let
$$b_v=-\frac1{1-|a_0|^2}\lt(\frac{\ll a_v,a_0\rl a_0}
{1+\sqrt{1-|a_0|^{2}}}+\sqrt{1-|a_0|^{2}}a_v\rt).$$ Using Lemma 1 to
the function $\phi$, by (2.3), we have $|b_v|^2v^v/|v|^{|v|}\le1$
and
$$\lt|\lt(\frac{\ll a_v,a_0\rl a_0}
{1+\sqrt{1-|a_0|^{2}}}+\sqrt{1-|a_0|^{2}}a_v\rt)\rt|^2
\le\frac{|v|^{|v|}}{v^v}(1-|a_0|^2)^2.$$ The same calculation as
in the proof of the above theorem gives
$$\lt|\lt(\frac{\ll a_v,a_0\rl a_0}{1+\sqrt{1-|a_0|^{2}}}
+\sqrt{1-|a_0|^{2}}a_v\rt)\rt|^2=\lt|\lt\ll a_v,\ a_0\rt\rl\rt|^2
+(1-|a_0|^2)\lt|a_v\rt|^2.$$ This shows (3.2).

Now, let the equality in (3.2) holds for some $v=(v_1,\cdots,v_n)$
with $v_j\neq0$ for $j=1,\cdots,n$. If $a_0=0$, then
$|a_v|^2v^v/|v|^{|v|}=1$, the equality in (2.3) holds and,
consequently, $f(z)=a_vz^v$. This shows (3.3). In the case
$a_0\ne0$, we have $|b_v|^2v^v/|v|^{|v|}=1$. Then, the same
reasoning shows $\phi(z)=b_vz^v$ and, consequently,
$$h(z)=\vp_{a_0}(b_vz^v)=\frac{a_0-\frac{\langle b_v,a_0\rangle}
{|a_0|^{2} }a_0z^v -\sqrt{1-|a_0|^{2}}\left(b_vz^v-\frac{\langle
b_v,a_0\rangle}{|a_0|^{2}}a_0z^v\right)}{1-\langle b_v,a_0\rangle
z^v}.\eqno(3.4)$$ Note that
$$\langle b_v,a_0\rangle=-\frac{\langle a_v,a_0\rangle}{1-|a_0|^{2}}.
\eqno(3.5)$$ Replacing $\langle b_v,a_0\rangle$ in (3.4) by (3.5),
by a straightforward calculation, we obtain
$$h(z)=a_0+\frac{a_v z^v}{1+\frac{\langle a_v, a_0\rangle
z^v}{1-|a_0|^2}}.$$

If $k=1$, $f(z)=h(z)$ and (3.3) is true. In the case $k\ge2$, we
have $f(z)=\vp_{a_0}(b_vz^v)+g(z)$ with
$$g(z)=\sum_{j=1}^{k-1}\sum_{m=0}^\infty\sum_{|\a|=mk+j} a_\a z^\a.$$
Let $0\le\te_1,\cdots,\te_n\le2\pi$ be fixed. For $\l\in\mb D$,
define
\begin{align*}
\psi(\l)&=f\left(e^{i\te_1}\sqrt{v_1}\l/\sqrt{|v|},\cdots,
e^{i\te_n}\sqrt{v_n}\l/ \sqrt{|v|}\right)\\
&=\vp_{a_0}\lt(be^{i(v_1\te_1+\cdots+v_n\te_n)}
\sqrt{\frac{v^v}{|v|^{|v|}}}\l^k\rt)\\
&\qquad\qquad+\sum_{j=1}^{k-1}\sum_{m=0}^\infty\left(\sum_{|\a|=mk+j}a_\a
e^{i(\a_1\te_1+\cdots+\a_n\te_n)}\cdot
\frac{\sqrt{v^\a}}{\sqrt{k^{mk+j}}}\right)\l^{mk+j}.
\end{align*}
Using Lemma 3 to $\psi$, we have
$$p_{m,j}(\te_1,\cdots,\te_n)=\sum_{|\a|=mk+j}a_\a
\sqrt{v^\a}e^{i(\a_1\te_1+\cdots+\a_n\te_n)}=0$$ for
$j=1,\cdots,k-1$ and $m=0,1,\cdots$. Note that the above equality
holds for arbitrary $\te_1,\cdots,\te_n$. Thus, for any multi-index
$\a'$ with $|\a'|=mk+j$, $1\le j\le k-1$, $m\ge0$, we have
$$a_{\a'}\sqrt{v^{\a'}}=\frac1{(2\pi)^n}
\int_0^{2\pi}\!\!\!\!\cdots\!\int_0^{2\pi}
e^{-i(\a'_1\te_1+\cdots+\a'_n\te_n)}
p_{m,j}(\te_1,\cdots,\te_n)d\te_1\cdots d\te_n=0.$$ It is proved
that $a_\a=0$ for any multi-index $\a$ with $|\a|=mk+j$, $1\le j\le
k-1$, $m\ge0$. Then we obtain $f(z)=h(z)$ and (3.3) is proved again.
The last conclusion of the theorem is easy to verify. The theorem is
proved. \ep

\noindent{\bf Remark 1.}\quad Define
$$f(z)=a_{1,0}z_1+a_{0,2}z_2^2=z_1+\frac13z^2_2,
\quad\mbox{for}\ \ z=(z_1,z_2)\in\mb B_2.$$ It is easy to verify
that $f\in \O_{2,1}$. Let $v=(1,0)$. We have $a_0=0, a_v=1$. $v,
a_0$, and $a_v$ satisfy the equality in (3.2), but $f(z)$ is not
expressed by (3.3). This example shows that the condition $v_j\neq0$
for $j=1,\cdots,n$ in the second part of the above theorem cannot be
omitted.

\bc Let $f(z)=\sum\limits_\a a_\a z^\a\in \O_{n,m}$. Then, for any
multi-index $v\neq0$,
$$|a_v|\leq\sqrt{\frac{|v|^{|v|}}{v^v}}\sqrt{1-|a_0|^2}$$
if $m\geq2$; and
$$|a_v|\le\sqrt{\frac{|v|^{|v|}}{v^v}}(1-|a_0|^2)$$
if $m=1$ or, more general, $\l_1a_0+\l_2 a_v=0$ with
$\l_1,\l_2\in\mb C$. \ec

\section{The Schwarz-Pick lemma of high order}

First we consider mappings from the unit disk into a unit ball $\mb
B_m$. The following theorem is the special case that $n=1$ of our
general Schwarz-Pick lemma of high order.

\bt Let $f\in\Omega_{1,m}$. Then,
$$|\langle f^{(k)}(z),
f(z)\rangle|^2+(1-|f(z)|^2)|f^{(k)}(z)|^2
\leq\lt[\frac{k!(1-|f(z)|^2)}{(1-|z|^{2})^{k}}(1+|z|)^{k-1}\rt]^2\eqno(4.1)$$
holds for $k\ge1$ and $z\in\mb D$.
\et

\bp Let $\xi\in\mb D$ and a positive integer $k$ be fixed. We
consider $g=f\circ\vp_\xi\in\O_{1,m}$, where
$$\vp_\xi(z)=\frac{\xi-z}{1-\ov \xi z}.$$
Let $g(z)=\overset{\infty}{\underset{l=0}{\sum}}c_lz^l$ with
$c_l=(c_{1,l},\cdots,c_{m,l})$ for $l=1,2,\cdots$. Then $c_0=f(\xi)$
and, by (3.2) for $n=1$,
$$|\langle c_l,c_0\rangle|^2+(1-|c_0|^2)|c_l|^2\leq(1-|c_0|^2)^2
\eqno(4.2)$$ holds for $l\ge1$.

It is easy to verify that
\begin{equation*}
\lt.\frac{d^l(\varphi_\xi(z)^j)}{dz^l}\rt|_{z=\xi}=
\begin{cases}
0, &l<j;\\
\frac{(-1)^j(\bar{\xi})^{l-j}}{(1-|\xi|^{2})^l}\frac{l!(l-1)!}
{(l-j)!(j-1)!}, &l\geq j.
\end{cases}
\end{equation*}
Let
$$A_j=\frac{(-1)^j\ov \xi^{k-j}}{(1-|\xi|^{2})^k}\frac{k!(k-1)!}
{(k-j)!(j-1)!}.$$ Since $f=g\circ\vp_\xi$, we have
$$f^{(k)}(\xi)=\sum_{j=1}^kc_jA_j,$$
and, using (4.2) and the Schwarz inequality,
\begin{equation*}
\begin{split}
|\langle f^{(k)}(\xi),f(\xi)\rangle|^2&
+(1-|f(\xi)|^2)|f^{(k)}(\xi)|^2\\
&=\lt|\sum_{j=1}^kA_j\langle c_j, c_0\rangle\rt|^2+
(1-|c_0|^2)\lt|\sum_{j=1}^kc_jA_j\rt|^2\\
&\leq\sum_{j=1}^k|A_j|\sum_{j=1}^k|A_j||\langle c_j,
c_0\rangle|^2+(1-|c_0|^2)\sum_{j=1}^k|A_j|\sum_{j=1}^k|A_j||c_j|^2\\
&=\sum_{j=1}^k|A_j|\sum_{j=1}^k|A_j|\lt(|\langle c_j,
c_0\rangle|^2+(1-|c_0|^2)|c_j|^2\rt)\\
&\le(1-|c_0|^2)^2\lt(\sum_{j=1}^k|A_j|\rt)^2.
\end{split}
\end{equation*}
On the other hand,
$$\sum_{j=1}^k|A_j|
=\frac{k!}{(1-|\xi|^{2})^k}\sum_{j=1}^k\frac{(k-1)!|\xi|^{k-j}}{(k-j)!(j-1)!}
=\frac{k!}{(1-|\xi|^{2})^k}(1+|\xi|)^{k-1}.$$ This shows (4.1). The
theorem is proved. \ep

\bc Let $f\in\Omega_{1,m}$. Then,  for $k\ge1$ and $z\in\mb D$,
$$|f^{(k)}(z)|\leq\frac{k!(1-|f(z)|^2)^{1/2}}{(1-|z|^{2})^{k}}(1+|z|)^{k-1};$$
and
$$|f^{(k)}(z)|\leq\frac{k!(1-|f(z)|^2)}{(1-|z|^{2})^{k}}(1+|z|)^{k-1}$$
if $\l_1f(z)+\l_2f^{(k)}(z)=0$ with $\l_1,\l_2\in\mb C$. \ec

\noindent{\bf Remark 2.}\quad In \cite{3}, the authors proved that
(1.1) is asymptotically sharp in the sense that for any two points
$z,w\in\mb D$, there exists a holomorphic function $f_{z,w}$ on $\mb
D$,  such that $f_{z,w}(z)=w$, $f_{z,w}(\mb D)\su\mb D$, and
$$\lim_{w\to\p \mb D}\frac{|f_{z,w}^{(k)}(z)|}{(1-|f_{z,w}(z)|^2)}
=\frac{k!(1+|z|)^{k-1}}{(1-|z|^{2})^{k}}$$ holds for any positive
integer $k$. In the same way, we can construct examples of mappings
to show (4.1) is also asymptotically sharp. For fixed points
$\xi\in\mb D\-\{0\}$, $\arg\xi=\te$, and $w\in\mb B_m\-\{0\}$, let
$b=-(1-|w|^2)/|w|$, and define
$$g_w(z)=\frac{w}{|w|}\frac{|w|-z}{1-|w|z}
=w(1+bz+b|w|z^2+b|w|^2z^3+\cdots),$$ and
$$f_w(z)=g_w\lt(-e^{-i\te}\frac{\xi-z}{1-\ov\xi z}\rt).$$
Then, $f_w(\xi)=w$, and
$$f_w^{(k)}(\xi)=\frac{-e^{-ki\te}k!(1-|w|^2)w}{|w|(1-|\xi|^2)^k}
\sum_{v=1}^k\frac{|w|^{v-1}(k-1)!|\xi|^{k-v}}{(v-1)!(k-v)!},$$
$$\frac{|\langle f_w^{(k)}(\xi),f_w(\xi)\rangle|^2
+(1-|f_w(\xi)|^2)|f_w^{(k)}(\xi)|^2}{(1-|f_w(\xi)|^2)^2}
=\frac{|f_w^{(k)}(\xi)|^2}{(1-|f_w(\xi)|^2)^2}$$
$$=\lt(\frac{k!}{(1-|\xi|^2)^k}\sum_{v=1}^k
\frac{|w|^{v-1}(k-1)!|\xi|^{k-v}}{(v-1)!(k-v)!}\rt)^2.$$ Thus,
$$\lim_{w\to\p\mb B_m}\frac{|\langle f_w^{(k)}(\xi),f_w(\xi)\rangle|^2
+(1-|f_w(\xi)|^2)|f_w^{(k)}(\xi)|^2}{(1-|f_w(\xi)|^2)^2}$$
$$=\lt(\frac{k!}{(1-|\xi|^2)^k}\sum_{v=1}^k
\frac{(k-1)!|\xi|^{k-v}}{(v-1)!(k-v)!}\rt)^2
=\lt(\frac{k!}{(1-|\xi|^2)^k}(1+|\xi|)^{k-1}\rt)^2.$$

Now, we are ready to prove our main result.

\bt Let $f\in\Omega_{n,m}$. Then,
\begin{align*}
H_{f(z)}(D_k&(f,z,\b),\ D_k(f,z,\b))\\
&\le k!^2\lt(1+\frac{|\ll\b,z\rl|}
{((1-|z|^2)|\b|^2+|\ll\b,z\rl|^2)^{1/2}}\rt)^{2(k-1)}(H_z(\b,\b))^k
\tag{4.3}
\end{align*}
holds for $k\ge1$, $\b\in\mb C^n\-\{0\}$ and $z\in\mb B_n$. Further,
in the case $n\le m$, the equality in $(4.3)$ holds for $k=1$, some
$z=\xi\in\mb B_n$, and any $\b\in\mb C^n$, i.e.,
$$H_{f(\xi)}(f'(\xi)\b,f'(\xi)\b)=H_\xi(\b,\b)\eqno(4.4)$$
holds for any $\b\in\mb C^n$, if and if
$F'(0)=\vp'_{f(\xi)}(f(\xi))f'(\xi)\vp'_\xi(0)$ satisfies
$\ov{F'(0)}^{T} F'(0)=I$, where $I$ is the identity matrix of
$n\times n$, and
$$f(z)=f(\xi)+\lt(\frac{1-\ll z,\xi\rl}{1-|\xi|^2}
+\frac{\ov{f(\xi)}^T f'(\xi)(z-\xi)}
{1-|f(\xi)|^2}\rt)^{-1}f'(\xi)(z-\xi).\eqno(4.5)$$

\et

\bp Let $k\ge1$, $\b=(\b_1,\cdots,\b_n)\in\mb C^n\-\{0\}$ and
$\xi=(\xi_1,\cdots,\xi_n)\in\mb B_n$ be given.  First assume that
$\b\in\p\mb B_n$. We consider the disk
$$\D=\{\z\in\mb C:|\xi+\z\b|^2=|\xi_1+\b_1\z|^2+\cdots+|\xi_n+\b_n\z|^2<1\}.$$
To make the equation of $\D$ clearer, let $U$ be a unitary matrix
such that $U\b=(1,0,\cdots,0)^T$. Denote $
U\xi=\eta=(\eta_1,\cdots,\eta_n)^T$. Here we identify a point in
$\mb C^n$ with a column matrix of $n\times 1$. Since
$$|\xi+\z\b|^2=|U(\xi+\z\b)|^2=|\eta_1+\z|^2+|\eta_2|^2+\cdots+|\eta_n|^2,$$
 we have
$$\D=\{\z\in\mb C:|\eta_1+\z|^2<1-|\eta_2|^2-\cdots-|\eta_n|^2\}.$$
Thus, if we set $\s=(1-|\eta_2|^2-\cdots-|\eta_n|^2)^{1/2}$,
$\g=\s\b$, and
$$\z=\s\o-\eta_1,
\quad z=L(\o)=\xi+\o\g-\eta_1\b,$$ $g(\o)=f(L(\o))$ is a holomorphic
mapping from $\mb D$ into $\mb B_m$.

Using (4.1) to the mapping $g$ and the point $\o=\o'=\eta_1/\s$, we
have
$$|\langle g^{(k)}(\o'),
g(\o')\rangle|^2+(1-|g(\o')|^2)|g^{(k)}(\o')|^2
\leq\lt[\frac{k!(1-|g(\o')|^2)}{(1-|\o'|^{2})^{k}}(1+|\o'|)^{k-1}\rt]^2.
$$ Note that $g(\o')=f(\xi)$, $|\eta|=|\xi|$,
$\eta_1=\ll\xi,\b\rl$ and
$$\s^2=1-|\eta|^2+|\eta_1|^2=1-|\xi|^2+|\ll\b,\xi\rl|^2,$$
$$|\o'|=\frac{|\ll\b,\xi\rl|}{(1-|\xi|^2+|\ll\b,\xi\rl|^2)^{1/2}},
\quad 1-|\o'|^{2}=\frac{1-|\xi|^2}{\s^2}.$$ By the chain rule,
$$g^{(k)}(\o')=\sum_{|\a|=k}\frac{k!}{\a!}\frac{\p^kf(\xi)}{\p z_1^{\a_1}\cdots\p z_n^{\a_n}
}\g^\a=\s^k \sum_{|\a|=k}\frac{k!}{\a!}\frac{\p^kf(\xi)}{\p
z_1^{\a_1}\cdots\p z_n^{\a_n} }\b^\a.$$ Thus,
$$\lt|\lt\ll\sum_{|\a|=k}\frac{k!}{\a!}\frac{\p^kf(\xi)}{\p z_1^{\a_1}\cdots\p z_n^{\a_n}
}\b^\a,\ f(\xi)\rt\rl\rt|^2
+(1-|f(\xi)|^2)\lt|\sum_{|\a|=k}\frac{k!}{\a!}\frac{\p^kf(\xi)}{\p
z_1^{\a_1}\cdots\p z_n^{\a_n} }\b^\a\rt|^2$$
$$\le k!^2(1-|f(\xi)|^2)^2\lt[\frac{1-|\xi|^2+|\ll\b,\xi\rl|^2}
{(1-|\xi|^2)^2}\rt]^k\lt(1+\frac{|\ll\b,\xi\rl|}
{(1-|\xi|^2+|\ll\b,\xi\rl|^2)^{1/2}}\rt)^{2(k-1)}.$$ (4.3) is proved
for $z=\xi$ and any $\b\in\p\mb B_n$. For a general $\b$, we may
consider $\b/|\b|$, since $(4.3)$ is homogeneous for $\b$. (4.3) is
proved completely.

Now assume that $n\le m$ and (4.4) holds for any $\b\in\mb C^n$.
Consider $F=\vp_{f(\xi)}\circ f\circ\vp_\xi$. By the invariance of
the Bergman metric, $H_0(F'(0)\b,F'(0)\b)=H_0(\b,\b)$, i.e.,
$|F'(0)\b|=|\b|$, holds for any $\b\in\mb C^n$. This shows that the
$m\times n$-matrix $F'(0)$ satisfies $\ov{F'(0)}^T F'(0)=I$, where
$I$ is the identity matrix of $n\times n$. Note that
$F'(0)=\vp'_{f(\xi)}(f(\xi))f'(\xi)\vp'_\xi(0)$. Thus, for $z\in\mb
B_n$, $F(z)=F'(0)z$ and
$$f(z)=\vp_{f(\xi)}\lt(\vp'_{f(\xi)}(f(\xi))f'(\xi)
\vp'_\xi(0)\vp_\xi(z)\rt).\eqno(4.6)$$ Using the formulas for
$\vp_a$ at the beginning of Section 2, we have
$$f'(\xi)\vp'_\xi(0)\vp_\xi(z)=\frac{(1-|\xi|^2)f'(\xi)(z-\xi)}
{1-\ll z,\xi\rl},$$
$$\vp'_{f(\xi)}(f(\xi))f'(\xi)\vp'_\xi(0)\vp_\xi(z)$$
$$=-\frac{1-|\xi|^2} {1-\ll
z,\xi\rl}\lt(\frac{(1-(1-|f(\xi)|^2)^{1/2})\ov{f(\xi)}^T
f'(\xi)(z-\xi)}{|f(\xi)|^2(1-|f(\xi)|^2)}f(\xi)
+\frac{f'(\xi)(z-\xi)} {(1-|f(\xi)|^2)^{1/2}}\rt),$$
$$\ll \vp'_{f(\xi)}(f(\xi))f'(\xi)\vp'_\xi(0)\vp_\xi(z),f(\xi)\rl=
-\frac{(1-|\xi|^2)\ov{f(\xi)}^T f'(\xi)(z-\xi)}{(1-\ll
z,\xi\rl)(1-|f(\xi)|^2)},\eqno(4.7)$$
$$P_{f(\xi)}(\vp'_{f(\xi)}(f(\xi))f'(\xi)\vp'_\xi(0)\vp_\xi(z))
=-\frac{(1-|\xi|^2)\ov{f(\xi)}^T f'(\xi)(z-\xi)}{|f(\xi)|^2(1-\ll
z,\xi\rl)(1-|f(\xi)|^2)}f(\xi).\eqno(4.8)$$ (4.5) follows from
(4.6), (4.7) and (4.8). Conversely, if
$A=\vp'_{f(\xi)}(f(\xi))f'(\xi)\vp'_\xi(0)$ satisfies $\ov{A}^T A=I$
and (4.5) holds, then
$$f(z)=\vp_{f(\xi)}\lt(\vp'_{f(\xi)}(f(\xi))f'(\xi)
\vp'_\xi(0)\vp_\xi(z)\rt)$$ and, by the invariance of the Bergman
metric,
$$H_{f(\xi)}(f'(\xi)\b,f'(\xi)\b)=H_{f(\xi)}(\vp'_{f(\xi)}(0)A\vp'_\xi(\xi)\b,
\vp'_{f(\xi)}(0)A\vp'_\xi(\xi)\b)$$
$$=H_0(A\vp'_\xi(\xi)\b,A\vp'_\xi(\xi)\b)=|A\vp'_\xi(\xi)\b|^2=|\vp'_\xi(\xi)\b|^2$$
$$=H_0(\vp'_\xi(\xi)\b,\vp'_\xi(\xi)\b)=H_\xi(\b,\b)$$
holds for any $\b\in\mb C^n$. The theorem is proved.
\ep

\section{Schwarz-Pick estimates for derivatives of any order}

On the basis of Theorem 4, we can deduce an estimate for partial
derivatives of arbitrary order of mappings in $\O_{n,m}$.

\bt Let $f\in\Omega_{n,m}$. Then,
$$\lt|\lt\ll\frac{\p^{|v|}f(z)}{\p z_1^{v_1}\cdots
\p z_n^{v_n} },\ f(z)\rt\rl\rt|^2+(1-|f(z)|^2)\lt|
\frac{\p^{|v|}f(z)}{\p z_1^{v_1}\cdots \p z_n^{v_n}}\rt|^2$$
$$\le\frac{|v|^{|v|}}{v^v}\lt[v!(1+|z|)^{|v|-1}\cdot
\frac{1-|f(z)|^2}{(1-|z|^2)^{|v|}}\rt]^2. \eqno(5.1)$$ holds for any
multi-index $v=(v_1,\cdots,v_n)\ne0$ and $z\in\mb B_n$. In
particular, if $f\in\O_{n,1}$, then (5.1) becomes
$$\lt|\frac{\p^{|v|}f(z)}{\p z_1^{v_1}\cdots \p z_n^{v_n}}\rt|
\le\sqrt{\frac{|v|^{|v|}}{v^v}}v!(1+|z|)^{|v|-1}\cdot
\frac{1-|f(z)|^2}{(1-|z|^2)^{|v|}}.\eqno(5.2)$$ \et

\bp Let $v=(v_1,\cdots,v_n)\ne0$ and $\xi\in\mb B_n$ be given, and
$k=|v|$. By (4.3),
$$\lt|\lt\ll\sum_{|\a|=k}\frac{k}{\a!}\frac{\p^kf(\xi)}{\p z_1^{\a_1}\cdots\p z_n^{\a_n}
}\frac{v^\a}{|v|},\ f(\xi)\rt\rl\rt|^2+(1-|f(\xi)|^2)\lt|\sum_
{|\a|=k}\frac{k}{\a!}\frac{\p^kf(\xi)}{\p z_1^{\a_1}\cdots\p
z_n^{\a_n}}\frac{v^\a}{|v|}\rt|^2\le A^2,$$ where
$$A=k!(1+|\xi|)^{k-1}\cdot\frac{1-|f(\xi)|^2}{(1-|\xi|^2)^k}.$$
Define
$$g(z)=\frac1A\lt\ll\sum_{|\a|=k}\frac{k!}{\a!}\frac{\p^kf(\xi)}
{\p z_1^{\a_1}\cdots\p z_n^{\a_n} }z^\a,\ f(\xi)\rt\rl
=\frac1A\sum_{|\a|=k}\lt\ll\frac{k!}{\a!}\frac{\p^kf(\xi)} {\p
z_1^{\a_1}\cdots\p z_n^{\a_n} },\ f(\xi)\rt\rl z^\a ,$$
$$h(z)=\frac1A(1-|f(\xi)|^2)^{1/2}\sum_{|\a|=k}\frac{k!}{\a!}\frac{\p^kf(\xi)}
{\p z_1^{\a_1}\cdots\p z_n^{\a_n}}z^\a,$$ and $\phi=(g,h)$. Using
(2.2) to $\phi$, which is a holomorphic mapping from $\mb B_n$ into
$\mb B_{2m}$ and satisfies $|\phi(z)|^2<1$ for $z\in\mb B_n$, we
have
$$\sum_{|\a|=k}\lt(\lt|\lt\ll\frac{k!}{\a!}\frac{\p^kf(\xi)}{\p z_1^{\a_1}\cdots
\p z_n^{\a_n} },\ f(\xi)\rt\rl\rt|^2+(1-|f(\xi)|^2)
\lt|\frac{k!}{\a!}\frac{\p^kf(\xi)}{\p z_1^{\a_1}\cdots \p
z_n^{\a_n} }\rt|^2\rt)\cdot\frac{v^\a}{|v|^{|\a|}}\le A^2.$$ In
particular,
$$\lt|\lt\ll\frac{k!}{v!}\frac{\p^kf(\xi)}{\p z_1^{v_1}\cdots
\p z_n^{v_n} },\ f(\xi)\rt\rl\rt|^2+(1-|f(\xi)|^2)
\lt|\frac{k!}{v!}\frac{\p^kf(\xi)}{\p z_1^{v_1}\cdots \p z_n^{v_n}
}\rt|^2\le\frac{|v|^{|v|}}{v^v}A^2.$$ This shows (5.1) and the
theorem is proved.
 \ep

\bt Let $f\in\Omega_{n,m}$. Then,
$$\lt|\lt\ll\frac{\p^{|v|}f(z)}{\p z_1^{v_1}\cdots
\p z_n^{v_n} },\ f(z)\rt\rl\rt|^2+(1-|f(z)|^2)\lt|
\frac{\p^{|v|}f(z)}{\p z_1^{v_1}\cdots \p z_n^{v_n}}\rt|^2$$
$$\le\frac{|v|^{|v|}}{v^v}\lt[v!\mu(z)\cdot
\frac{1-|f(z)|^2}{(1-|z|^2)^{(v_1+|v|)/2}}\rt]^2 \eqno(5.3)$$ holds
for any multi-index $v=(v_1,\cdots,v_n)\ne0$ and
$z=(z_1,0,\cdots,0)\in\mb B_n$, where $\mu(z)=(1+|z|)^{|v|-1}$ if
$v_1=|v|$, and $\mu(z)$ is the sum of terms $c_j|z|^j$ with $j\le
v_1$ in $(1+|z|)^{|v|-1}$. In particular, if $f\in\O_{n,1}$, then
(5.3) becomes
$$\lt|\frac{\p^{|v|}f(z)}{\p z_1^{v_1}\cdots \p z_n^{v_n}}\rt|
\le\sqrt{\frac{|v|^{|v|}}{v^v}}v!\mu(z)\cdot
\frac{1-|f(z)|^2}{(1-|z|^2)^{(v_1+|v|)/2}}.\eqno(5.4)$$ \et

\bp Let $v=(v_1,\cdots,v_n)\ne0$ and $\xi=(\xi_1,\cdots,\xi_n)\in\mb
B_m$ be given. $v_1=|v|$, (5.3) follows from (5.1). Now assume that
$v_1<|v|$. Let $k=|v|$ and
$$g(z)=f(\vp_\xi(z))=\sum_{\a}c_\a z^\a.$$
Then, $c_0=f(\xi)$ and, by (3.2),
$$\lt|\lt\ll c_\a,\ c_0\rt\rl\rt|^2
+(1-|c_0|^2)\lt|c_\a\rt|^2\le\frac{|\a|^{|\a|}}{\a^a}(1-|c_0|^2)^2\eqno(5.5)$$
holds for any multi-index $\a\ne0$. Thus, we have
$$f(z)=g(\vp_\xi(z))=\sum_{\a}c_\a\vp_\xi(z)^\a,$$
where
$$\vp_\xi(z)=\lt(\frac{\xi_1-z_1}{1-\ov\xi_1z_1},-\frac{(1-|\xi|^2)^{1/2}z_2}
{1-\ov\xi_1z_1},\cdots,-\frac{(1-|\xi|^2)^{1/2}z_n}{1-\ov\xi_1z_1}\rt).$$
For a multi-index $\a$, denote $\a=(\a_1,\a')$ with
$\a'=(\a_2,\cdots,\a_n)$. Then, it is easy to see that
\begin{equation*}
\lt.\frac{\p^k(\vp_\xi(z)^\a)}{\p z_1^{v_1}\cdots \p
z_n^{v_n}}\rt|_{z=\xi}=
\begin{cases}
0, &\a'\ne v';\\
0,&\a'\ne v',\ \a_1>v_1;\\
\frac{(-1)^{\a_1+|v'|}(\bar{\xi})^{v_1-\a_1}}{(1-|\xi|^{2})^{(v_1+|v|)/2}}\frac{v!(k-1)!}
{(v_1-\a_1)!(\a_1-1+|v'|)!}, &\a'=v',\ 0\le\a_1\le v_1.
\end{cases}
\end{equation*}
Thus, letting
$$A_j=\frac{(-1)^{j+|v'|}(\bar{\xi})^{v_1-j}}
{(1-|\xi|^{2})^{(v_1+|v|)/2}}\frac{v!(k-1)!}{(v_1-j)!(j-1+|v'|)!},$$
we have
$$\lt.\frac{\p^kf(z)}{\p z_1^{v_1}\cdots \p z_n^{v_n}}\rt|_{z=\xi}
=\sum_{j=0}^{v_1}A_jc_{j,v'},$$ and, by the Schwarz inequality and
(5.5),
$$\lt|\lt\ll\lt.\frac{\p^{|v|}f(z)}{\p z_1^{v_1}\cdots
\p z_n^{v_n}}\rt|_{z=\xi},\
f(\xi)\rt\rl\rt|^2+(1-|f(\xi)|^2)\lt|\lt. \frac{\p^{|v|}f(z)}{\p
z_1^{v_1}\cdots \p z_n^{v_n}}\rt|_{z=\xi}\rt|^2$$
$$=\lt|\sum_{j=0}^{v_1}A_j\ll c_{j,v'},c_0\rl\rt|^2+(1-|c_0|^2)
\lt|\sum_{j=0}^{v_1}A_jc_{j,v'}\rt|^2$$
$$\le\sum_{j=0}^{v_1}|A_j|\sum_{j=0}^{v_1}|A_j||\ll c_{j,v'},c_0\rl|^2
+(1-|c_0|^2)\sum_{j=0}^{v_1}|A_j|\sum_{j=0}^{v_1}|A_j||c_{j,v'}|^2$$
$$\le\frac{|v|^{|v|}}{v^v}(1-|c_0|^2)^2\lt(\sum_{j=0}^{v_1}|A_j|\rt)^2.$$
Here, we use the obvious inequality
$\frac{|\a|^{|\a|}}{\a^\a}\le\frac{|v|^{|v|}}{v^v}$ if $\a_j\le v_j$
for $j=1,\cdots,n$. Note that
$$\sum_{j=0}^{v_1}|A_j|=\frac{v!}{(1-|\xi|^{2})^{(v_1+|v|)/2}}
\sum_{j=0}^{v_1}\frac{(k-1)!|\xi|^{v_1-j}}{(v_1-j)!(j-1+|v'|)!}$$
$$=\frac{v!}{(1-|\xi|^{2})^{(v_1+|v|)/2}}\sum_{l=0}^{v_1}\frac{(k-1)!
|\xi|^l}{l!(k-1-l)!}$$ (5.3) is proved. (5.4) follows from (5.3)
directly and the proof is complete. \ep

\noindent{\bf Remark 3.}\quad Let a multi-index
$v=(v_1,\cdots,v_n)$, $\xi=(\xi_1,0,\cdots,0)\in\mb B_n$ and
$w\in\mb C$ be fixed. Define
$$g(z)=\frac{w-\sqrt{|v|^{|v|}/v^v}z^v}{1-\ov w\sqrt{|v|^{|v|}/v^v}z^v}\quad\mbox{for}
\ \ z\in\mb B_n,$$ and $f=g\circ\vp_\xi$. Then, $g\in\O_{n,1}$,
$f\in\O_{n,1}$, $f(\xi)=w$, and
$$g(z)=w-(1-|w|^2)\sqrt{|v|^{|v|}/v^v}z^v-\ov w(1-|w|^2)(|v|^{|v|}/v^v)z^{2v}
-\cdots,$$
$$\lt.\frac{\p^{|v|}f(z)}{\p z_1^{v_1}\cdots \p z_n^{v_n}}\rt|_{z=\xi}
=-\sqrt{\frac{|v|^{|v|}}{v^v}}(1-|w|^2)\cdot\lt.\frac{\p^{|v|}(\vp_\xi(z))^v}{\p
z_1^{v_1}\cdots\p z_n^{v_n}}\rt|_{z=\xi}$$
$$=(-1)^{|v|+1}\sqrt{\frac{|v|^{|v|}}{v^v}}v!\cdot\frac{1-|f(\xi)|^2}
{(1-|\xi|^2)^{(v_1+|v|)/2}}.$$ This shows that the estimate (5.4) is
precise up to a constant less than $2^{|v|-1}$.
\medskip

\noindent{\bf Remark 4.}\quad  If $v_1=|v|=k$, (5.4) becomes
$$\lt|\frac{\p^kf(z)}{\p z_1^k}\rt|
\le k!(1+|z|)^{k-1}\cdot \frac{1-|f(z)|^2}{(1-|z|^2)^k}.\eqno(5.5)$$
(5.5) is also a consequence of (5.2). For given
$\xi=(\xi_1,0,\cdots,0)\in\mb B_n\-\{0\}$ and $w\in\mb C\-\{0\}$,
let $\te=\arg\xi-\arg w$ and defined
$$g(z)=\frac{w+e^{-i\te}z_1}{1+\ov we^{-i\te}z_1}\quad\mbox{for}\ \ z=(z_1,\cdots,z_n)\in\mb B_n,$$
and $f=g\circ\vp_\xi$. Then, $g\in\O_{n,1}$, $f\in\O_{n,1}$,
$f(\xi)=w$,  and
$$g(z)=w+(1-|w|^2)e^{-i\te}z_1-\ov w(1-|w|^2)e^{-2i\te}z_1^2
+\ov w^2(1-|w|^2)e^{-3i\te}z_1^3+\cdots.$$ Thus, for any positive
integer $k$, we have
$$\lt.\frac{\p^kf(z)}{\p z_1^k}\rt|_{z=\xi}=(1-|w|^2)\sum_{j=1}^k\ov w^{j-1}
e^{-ji\te}\lt.\frac{d^k}{dz_1^k}\lt(\frac{(\xi_1-z_1)^j}{(1-\ov\xi_1z_1)^j}\rt)
\rt|_{z_1=\xi_1}$$
$$=-\frac{k!(1-|w|^2)}{(1-|\xi|^2)^k}\sum_{j=1}^k\frac{(k-1)!e^{-ji\te}\ov w^{j-1}
\ov\xi^{k-j}}{(j-1)!(k-j)!}$$
$$=-\frac{k!(1-|w|^2)|\xi|^kw}{(1-|\xi|^2)^k|w|\xi^k}\sum_{j=1}^k\frac{(k-1)|w|^{j-1}|\xi|^{k-j}}
{(j-1)!(k-j)!}$$
$$=-\frac{k!(1-|w|^2)|\xi|^kw}{(1-|\xi|^2)^k|w|\xi^k}(|w|+|\xi|)^{k-1}$$ and, consequently,
$$\lim_{w\to\p\mb D}\frac{\lt|\lt.\p^kf(z)/\p z_1^k\rt|_{z=\xi}\rt|}
{1-|f(\xi)|^2}= \frac{k!(1+|\xi|)^{k-1}}{(1-|\xi|^2)^k}.$$ This
shows that (5.5) is asymptotically sharp.

\end{document}